\documentclass[11pt]{article}
\usepackage[french]{babel}
\usepackage[latin1]{inputenc}
\usepackage{amsmath}
\usepackage{amscd}
\usepackage{amstext}
\usepackage{amsbsy}
\usepackage{amsopn}
\usepackage{amsthm}
\usepackage{amsxtra}
\usepackage{amssymb}
\usepackage{hyperref}
\usepackage[all]{xy}
\newcommand{\be}{\begin{enumerate}}
\newcommand{\ee}{\end{enumerate}}
\newtheorem{cor}[subsubsection]{Corollaire}

\newtheorem{sousprop}[paragraph]{Proposition}

\newtheorem{prop}[subsubsection]{Proposition}

\newtheorem{sousdef}[paragraph]{D\'efinition}

\newtheorem{surthm}[subsection]{Théorème}

\newtheorem{thm}[subsubsection]{Th\'eor\`eme}
\newcommand{\inft}{{\scriptstyle\infty}}

\newcommand{\what}{\widehat}
\newcommand{\ot}{\otimes}

\newcommand{\ti}{\tilde}

\newcommand{\rig}{\rightarrow}
\newcommand{\trig}{\twoheadrightarrow}

\newcommand{\hrig}{\hookrightarrow}
\newcommand{\sta}{\stackrel}

\newcommand{\gr}{{\rm gr}}
\newcommand{\pg}{_{\bullet}}

\newcommand{\crofrac}[2]{\genfrac{\langle}{\rangle}{0pt}{}{#1}{#2}}

\newcommand{\acfrac}[2]{\genfrac{\{}{\}}{0pt}{}{#1}{#2}}

\newcommand{\der}{\partial}
\newcommand{\la}{\langle}
\newcommand{\ra}{\rangle}

\newcommand{\spf}{{\rm Spf}\,}
\newcommand{\Qr}{{\bf Q}}
\newcommand{\Ne}{{\bf N}}
\newcommand{\Ze}{{\bf Z}}
\newcommand{\Cc}{{\bf C}}

\newcommand{\bP}{{\bf P}}
\newcommand{\Gm}{{\bf G}_m}

\newcommand{\al}{\alpha}
\newcommand{\varep}{\varepsilon}
\newcommand{\lam}{\lambda}

\newcommand{\Ga}{\Gamma}

\newcommand{\Lam}{\Lambda}
\renewcommand{\AA}{{\cal A}}

\newcommand{\CC}{{\cal C}}
\newcommand{\DD}{{\cal D}}
\newcommand{\EE}{{\cal E}}
\newcommand{\FF}{{\cal F}}
\newcommand{\GG}{{\cal G}}
\newcommand{\II}{{\cal I}}
\newcommand{\HH}{{\cal H}}
\newcommand{\JJ}{{\cal J}}

\newcommand{\LL}{{\cal L}}
\newcommand{\MM}{{\cal M}}
\newcommand{\NN}{{\cal N}}
\newcommand{\OO}{{\cal O}}

\renewcommand{\SS}{{\cal S}}
\newcommand{\TT}{{\cal T}}

\newcommand{\XX}{{\cal X}}

\newcommand{\ZZ}{{\cal Z}}

\newcommand{\uk}{\underline{k}}
\newcommand{\ul}{\underline{l}}

\newcommand{\uy}{\underline{y}}

\newcommand{\uder}{\underline{\der}}

\newcommand{\Dm}{\DD^{(m)}}
\newcommand{\Dcm}{\what{\DD}^{(m)}}
\newcommand{\Ddag}{\DD^{\dagger}}
\newcommand{\spec}{{\rm spec}\,}

\newcommand{\bS}{{\mathbf S}}
\newcommand{\bSm}{{\mathbf S}^{(m)}}

\newcommand{\bGamm}{{\mathbf \Gamma}_{(m)}}
\begin{document}
%
\title{Un théorème de Beilinson-Bernstein pour les $\DD$-modules arithmétiques}
\author{C.~Noot-Huyghe\,\footnote{This work has been supported by the research network
Arithmetic Algebraic Geometry of the European Community (Programme FP6, contrat
MRTN-CT2003-504917)}}
\maketitle
\begin{abstract}
In the 80's, Brylinski-Kashiwara and Beilinson-Bernstein
proved that flag varieties over {\bf C} are $\DD$-affine.
We show an analogous of this theorem for
arithmetic $\DD$-modules (in the sense of Berthelot) over flag varieties
defined over the formal spectrum of a complete discrete valuation ring of 
inequal characteristics.

Dans les années 80, Brylinski-Kashiwara et Beilinson-Bernstein 
ont démontré que les variétés de drapeaux sur {\bf C} sont 
$\DD$-affines. Nous montrons un analogue de ce théorème pour les 
$\DD$-modules arithmétiques (au sens de Berthelot) sur les variétés de 
drapeaux sur le spectre formel d'un anneau de valuation discrète complet
d'inégales caractéristiques.
\end{abstract}
\tableofcontents
\section*{Introduction}
Soit $V$ un anneau de valuation discrète, d'inégales caractéristiques $(0,p)$.
On considère ici les deux situations suivantes: 
\be \item[1-] $S=\spec V$, le spectre de $V$, $X$ est un $S$-schéma noetherien,
\item[2-]  $\SS=\spf V$, le spectre formel de $V$, $\XX$ un schéma formel noetherien sur 
$\SS$. \ee
Soit $\AA$ un faisceau cohérent de $\OO_X$-modules (resp. un faisceau cohérent de
$\OO_{\XX}$-modules). Un $\AA$-module sur le schéma $X$ sera dit quasi-cohérent s'il est un $\OO_X$-module quasi-cohérent. On
dit que $X$ (resp. $\XX$) est $\AA$-affine si les deux propriétés suivantes sont 
vérifiées:
\be \item[(i)] Pour tout $\AA$-module quasi-cohérent $\MM$ sur $X$ (resp. tout $\AA$-module 
cohérent sur $\XX$) et tout $n\geq 1$ on a les égalités $H^n(X,\MM)=0$ (resp.
$H^n(\XX,\MM)=0$). 
\item[(ii)] Le foncteur $\Gamma$ établit une équivalence de catégories entre la catégorie des 
$\AA$-modules quasi-cohérents (resp. des $\AA$-modules cohérents) et la catégorie des 
$\Ga(X,\AA)$-modules (resp. des $\Ga(\XX,\AA)$-modules de type fini).
\ee

Un énoncé important de la théorie des groupes est le théorème de Beilinson-Bernstein: 
soit $G$ un groupe semi-simple sur $\Cc$, $X$ la variété de drapeaux de $G$, $\DD_X$ le 
faisceau des opérateurs différentiels sur $\XX$, alors $X$ est $\DD_X$-affine. On se
propose de donner ici un analogue arithmétique de cet énoncé, dans la situation qui suit. 
Soit  $G$ un groupe semi-simple sur $S$, $\rho$ la demi-somme des racines
positives de $G$, $P$ un sous-groupe parabolique de $G$, 
$X=G/P$, qu'on suppose défini sur $S$, $\XX$ le schéma 
formel obtenu en complétant $X$ le long de la fibre spéciale de $S$. Ce schéma est lisse et on peut s'intéresser 
au faisceau des opérateurs différentiels arithmétiques sur $\XX$ construit par Berthelot,
que nous noterons $\Ddag_{\XX,\Qr}$. On montre alors que $\XX$ est $\Ddag_{\XX,\Qr}$-affine. Plus 
généralement, $\XX$ est $\Ddag_{\XX,\Qr}(\lam)$-affine pour tout poids $\lam$ tel que 
$\lam+\rho$ est dominant et régulier, le faisceau $\Ddag_{\XX,\Qr}(\lam)$ désignant le 
faisceau des opérateurs différentiels arithmétiques à valeurs dans $\OO_{\XX}(\lam)$.

En caractéristique $0$, pour le faisceau $\DD(\lam)$ tel que $\lam+\rho$ est régulier, le 
résultat est démontré indépendamment par Beilinson-Bernstein (\cite{BeBe}) et 
par Brylinski-Kashiwara (\cite{Bry-Kash}) et joue un rôle essentiel dans la démonstration
de la conjecture de multiplicité de Kazhdan-Lusztig (\cite{Kazhdan-Lusztig1}). 

En caractéristique $p>0$, Haastert a montré que cet énoncé d'affinité était 
vérifié pour les espaces projectifs, ainsi que pour la variété de drapeaux de 
$SL_3$ (\cite{Haa}). En revanche, Kashiwara-Lauritzen ont donné un contre-exemple à cet énoncé, pour le 
faisceau usuel $\DD$ (\cite{Kash-Laur}) et pour la grassmanienne des sous-espaces 
vectoriels de dimension $2$ d'un espace de dimension $5$. Enfin, 
Bezrukavnikov, Mirkovic, Rumynin ont montré (3.2 de~\cite{Bezru1}) un analogue de ce résultat 
d'affinité en passant à la 
 catégorie dérivée bornée des $\DD^{(0)}$-modules cohérents sur $X$ (i.e. les opérateurs 
différentiels sans puissances divisées) et sous la condition que 
$p$ soit strictement plus grand que le nombre de Coxeter de $G$. 

En caractéristique mixte, le résultat a été montré pour les espaces projectifs
(\cite{Hu1}). Dans ce cas, on utilise de façon cruciale le fait que 
le faisceau tangent est très ample, ce qui caractérise l'espace projectif. Le point clé pour les variétés de drapeaux est 
que la catégorie des $\DD^{\dagger}_{\XX,\Qr}$-modules cohérents est engendrée par les
modules induits (i.e. du type $\DD^{\dagger}_{\XX,\Qr}\ot_{\OO_{\XX}}\EE$ où $\EE$ est 
un $\OO_{\XX}$-module cohérent). On utilise cette propriété pour montrer que 
si le résultat de $\DD$-affinité est vrai algébriquement, pour le 
faisceau $\DD_{X_K}$, alors il est vrai pour le faisceau 
$\DD^{\dagger}_{\XX,\Qr}$ sur le schéma formel $\XX$ (théorème~\ref{subsection-thm-int}).

Nous n'aborderons pas ici l'aspect localisation de $Lie(G)$-modules (ou plutôt des modules
sur la complétion faible de $Lie(G)$), qui est 
bien entendu sous-jacent et fera l'objet d'un article ultérieur.
\section{Notations-Rappels}
\label{section-rappels}
\subsection{Notations}
Dans toute la suite, on note $K$ le corps des
fractions de $V$, $\pi$ une uniformisante et $k$ le corps résiduel de $V$. 
Soit $G$ un groupe semi-simple sur $S$, 
$P$ un sous-groupe parabolique de $G$, $X=G/P$, $\XX$ le schéma formel associé par
complétion à $X$.

 D'une façon générale, si $Z$ est un $S$-schéma, la
lettre cursive $\ZZ$ désignera le schéma formel obtenu en complétant $Z$ le long de
l'idéal $\pi$, $Z_k$ la fibre spéciale $Z_k=\spec k \times_S Z$ et 
$Z_K$ la fibre générique $Z_K=\spec K \times_S Z$. On posera 
$i_Z$ l'immersion fermée $Z_k \hrig Z$ et $j_Z$ l'immersion ouverte 
$Z_K\rig Z$ (on omettra éventuellement le $Z$ dans les notations quand le contexte sera
clair). On notera aussi 
$$Z_i=Z\times_{S}\spec (S/\pi^{i+1}S).$$ Pour un faisceau $\EE$ sur 
un $S$-schéma $Z$, on notera $\EE_k=i^*\EE$ et $\EE_K=j^*\EE$.
\subsection{Coefficients $p$-adiques}
Fixons un entier $m$. Si $k_i\in\Ne$, on introduit 
$q_{k_i}$ le quotient de la division euclidienne de $k_i$ par $p^m$ et 
pour un multi-indice $\uk=(k_1,\ldots,k_N)$ on définit 
$$q_{\uk}!=\prod_{i=1}^N q_{k_i}!.$$ 

Pour $k\leq l\in\Ne$, on pose
$$\acfrac{l}{k}=\frac{q_k!}{q_k!q_{l-k}!},$$
$$\crofrac{l}{k}=\binom{l}{k}\acfrac{l}{k}^{-1}\in {\bf Z}_{(p)},$$
et pour des multi-indices $\uk$, $\ul \in \Ne^N$, tels que $\uk\leq \ul$ 
(i.e. $k_i\leq l_i$ pour tout $1\leq i \leq N$), 
$$\crofrac{\ul}{\uk}=\prod_{i=1}^N\crofrac{l_i}{k_i}.$$
On définit de façon analogue les coefficients 
$\acfrac{\ul}{\uk}$ et $\binom{\ul}{\uk}$. 

\vspace{+3mm}
Décrivons maintenant les différents
faisceaux d'opérateurs différentiels intervenant dans cette situation.
%
\subsection{Opérateurs différentiels arithmétiques}
Dans cette partie, $X$ est un $S$-schéma formel lisse et $\XX$ est son complété formel
lelong de l'idéal engendré par $\pi$. 
On décrit les différents faisceaux d'opérateurs différentiels en coordonnées locales. Nous renvoyons à 
\cite{EGA4} et à \cite{Be1} pour une définition intrinsèque de ces faisceaux. Soit 
$U$ un ouvert affine lisse de $X$, $x_1,\ldots,x_N$ une famille de coordonnées locales sur 
$X$, $dx_1,\ldots,dx_N$ une base de $\Omega_X^1(U)$, $\der_1,\ldots,\der_N$ la base duale 
de $\TT_X(U)$. Si $k_i\in\Ne$, on note $\der_i^{[k_i]}=\der_i/k_i!$ et pour un
multi-indice $\uder^{[\uk]}=\prod_{i=1}^N\der_i^{[k_i]}$. Alors on a la description 
suivante (\cite{EGA4})
$$\DD_{X}(U)=\left\{\sum_{finies}a_{\uk} \uder^{[\uk]}\,|\,a_{\uk}\in\OO_X(U)\right\}.$$
Donnons maintenant une description des faisceaux d'opérateurs différentiels construits par 
P.~Berthelot. 

Soit $m\in\Ne$. P.~Berthelot introduit les faisceaux $\Dm_{X}$, ainsi que $\Dcm_{\XX}$,
leur complété $p$-adique sur $\XX$.  Notons
$$ \uder^{{\la \uk \ra}_{(m)}}=q_{\uk}!\uder^{[\uk]}.$$
On a alors les descriptions suivantes
$$\DD^{(0)}_{X}(U)=\left\{\sum_{finies}a_{\uk} \uder^{\uk}\,|\,a_{\uk}\in\OO_X(U)\right\},$$
$$\DD^{(m)}_{X}(U)=\left\{\sum_{finies}a_{\uk} \uder^{\la \uk
\ra_{(m)}}\,|\,a_{\uk}\in\OO_X(U)\right\}.$$
Pour $m=+\inft$, on retrouve le faisceau usuel $\DD_X$. Les faisceaux $\Dm_{X}$ forment un 
système inductif, ainsi que leurs complétés $p$-adiques $\Dcm_{\XX}$. On posera 
$$\Ddag_{\XX,\Qr}=  \varinjlim_m \Dcm_{\XX,\Qr}.$$

Les faisceaux $\Dm_{X}$ sont à sections noethériennes sur les ouverts affines, c'est
donc aussi le cas des faisceaux $\Dcm_{\XX,\Qr}$, qui sont cohérents. On en déduit, via 
un théorème de platitude que le faisceau $\Ddag_{\XX,\Qr}$ est cohérent. 

Plus précisément, 
la structure de l'algèbre graduée de $\Dm_{X}$ est décrite en 1.3.7.3 de ~\cite{Hu1} en termes 
d'algèbre symétrique de niveau $m$ du faisceau tangent. 

Rappelons comment est construite l'algèbre symétrique de niveau $m$ d'un $\OO_X$-module 
localement libre $\EE$ (section 1 de \cite{Hu1}). 

\subsection{Algèbres symétriques de niveau $m$.}
\label{subsection-sym_m}
Pour les définitions relatives aux 
$m$-PD-structures, on se reportera à \cite{Be1}. Soit $\EE$ un $\OO_X$-module localement
libre. 
Le faisceau d'algèbres symétriques 
$\bS(\EE^{\vee})$ est gradué et muni de l'idéal d'augmentation $I(\EE^{\vee})=\bigoplus_{n\geq
1}S_n(\EE^{\vee})$. Par définition,
$\bGamm(\EE^{\vee})$ est la $m$-PD-enveloppe du couple $(\bS(\EE^{\vee}),I(\EE^{\vee}))$. Ce faisceau 
d'algèbres est muni d'un $m$-PD-idéal $\overline{\II}$, définissant une $m$-PD-filtration et 
on définit
$$\Ga_{(m)}^n(\EE^{\vee})=\bGamm(\EE^{\vee})/\overline{\II}^{\left\{n+1\right\}}.$$
On pose enfin 
$$\bS^{(m)}(\EE)=\bigcup_n \HH om_{\OO_X}(\Ga_{(m)}^n(\EE^{\vee}),\OO_X),$$
qui est (1.3.3 de \cite{Hu1}) un faisceau de $\OO_X$-algèbres commutatives graduées
par $$\bS^{(m)}(\EE)=\bigoplus_{n\in\Ne}S_n^{(m)}(\EE), \  \rm{o\grave{u}}\  
S_n^{(m)}(\EE)=\HH om_{\OO_X}(\overline{\II}^{\left\{n\right\}}/\overline{\II}^{\left\{n+1\right\}},\OO_X).$$
Les modules $S_n^{(m)}(\EE)$ sont localement libres de rang fini et le faisceau 
$\bS^{(m)}(\EE)$ est un faisceau d'algèbres localement noetheriennes. Ces constructions 
définissent des foncteurs contravariants $\bGamm$ et covariants $\bSm$ de la catégorie
des $\OO_X$-modules localement libres vers la catégorie des faisceaux d'algèbres
commutatives graduées. 

Si $y_1,\ldots y_N$ sont une base locale de $\EE$ sur un ouvert $U$ de $X$, le 
faisceau $Sm_n(\EE)$ admet pour base sur $U$ des éléments 
$$\uy^{\la \uk \ra}=y_1^{\la k_1\ra}y_2^{\la k_2\ra}\cdots y_N^{\la k_N\ra} \quad {\rm
tels}\quad {\rm que}\quad |\uk|= n.$$ De plus ces éléments vérifient
$$\uy^{\la \uk \ra}\cdot \uy^{\la \ul \ra}= \crofrac{\uk+\ul}{\uk}\uy^{\la \uk+\ul\ra}.$$

Les algèbres $\bS^{(m)}(\EE)$ vérifient les propriétés usuelles des algèbres 
symétriques. Nous aurons besoin de la propriété de dévissage suivante (1.3.9 de \cite{Hu1}).
Soit $0 \rig \EE \rig \FF \rig \GG \rig 0$ une suite exacte de
$\OO_X$-modules localement libres. On pose, pour $0\leq l \leq k$,
$$\Lam_{k}^l=\sum_{i\geq l}Im(\bS_i^{(m)}(\EE)\ot_{\OO_X}\bS_{k-i}^{(m)})(\FF)\rig
\bS_k^{(m)}(\FF)).$$
Les modules $\Lam_k^l$ forment une filtration décroissante de $S_k^{(m)}(\FF)$. Les
modules $\Lam_k^0$ et $\Lam_k^k$ sont isomorphes respectivement à 
$S_k^{(m)}(\FF)$ et $S_k^{(m)}(\EE)$.
\begin{prop} \label{subsubsection-devis}
Pour tout $l\leq k$, il existe des suites exactes de $\OO_X$-modules 
$$0\rig \Lam_k^{l+1}\rig \Lam_k^l\rig S^{(m)}_l(\EE)\ot_{\OO_X}S^{(m)}_{k-l}(\GG)\rig 0.$$
\end{prop} 
 L'intérêt de cette construction pour nous est le résultat suivant (1.3.7.3 de
\cite{Hu1}). 
\begin{prop}\label{subsubsection-isom_alg_grad}Il existe un isomorphisme canonique de faisceaux de $\OO_X$-algèbres graduées
$$gr_{\pg}\Dm_X \simeq \bS^{(m)}(\TT_X).$$
\end{prop} 
\subsection{Opérateurs différentiels à valeurs dans un faisceau inversible}
\label{subsection-diff_finv}
Soit $\LL$ un faisceau inversible sur $X$ (resp. un faisceau inversible sur $\XX$). On note
$\sharp$ l'un des symboles $(m)$ ou $\dagger$
$$\DD^{(m)}_{X}(\LL)=\LL \ot_{\OO_X}\DD^{(m)}_{X} \ot_{\OO_X}\LL^{-1},\ 
 \DD^{\sharp}_{\XX}(\LL)=\LL \ot_{\OO_{\XX}}\DD^{\sharp}_{\XX}
\ot_{\OO_{\XX}}\LL^{-1}\quad ({\rm resp.}\quad \what{\DD}^{(m)}_{\XX}(\LL)).$$
C'est le faisceau des opérateurs différentiels à valeurs dans $\LL$.
Si $\LL$ est un faisceau inversible sur $X$, on notera toujours $\LL$ le
$\OO_{\XX}$-module localement libre de rang $1$ obtenu en complétant $\LL$ le long 
de $\pi$.

\vspace{+3mm}
Donnons maintenant quelques considérations sur les liens entre système de 
racines de $Lie(G)$ et donnée de racines de $G$.
\subsection{Comparaison de systèmes de racines.}
Dans cette partie, on suppose que $G$ est semi-simple déployé. On pourra remplacer 
$V$ par $\Ze$ si $G$ est défini sur $\Ze$ et un tore maximal est déployé sur $\Ze$. 
On peut alors introduire la donnée de racines de ce groupe algébrique. On peut aussi
considérer le système de racines de l'agèbre de Lie de $G$. On explique ici comment 
identifier ces données (après avoir tensorisé par $K$). Dans le cas complexe, ces résultats sont bien connus. Faute de référence 
dans notre cas, nous expliquons comment procéder. 
\subsubsection{Donnée de racines d'un groupe algébrique et algèbre de distributions}
On note $1_G$ l'élément neutre de $G$, 
$1_G$: $\spec V\hrig G$, $\varep$ l'application correspondante: $V[G]\rig V[G]$, 
$i_G$: $G\rig G$ l'application de passage à l'inverse et $\sigma_G$: $V[G]\rig V[G]$
l'application correspondante. 
Soit $T$ un tore maximal déployé fixé de $G$.
La présentation de Jantzen (\cite{Jantzen}) est particulièrement bien adaptée à notre cadre. 
L'algèbre de groupe $V[G]$ se décompose $V[G]=V\bigoplus I_1$ où 
$I_1$ est l'idéal d'augmentation de $G$, c'est-à-dire le noyau du morphisme 
$\varep$. On pose alors 
$$Lie(G)=Hom_V(I_1/I_1^2,V).$$ Si $G$ est lisse, on introduit 
$\TT_G$ le faisceau tangent du groupe $G$. La suite exacte des faisceaux 
de formes différentielles appliquée à l'immersion fermée $1_G$ 
(voir par exemple chapitre II prop. 7 de \cite{Neron_models}) donne 
que $J_1/J_1^2\simeq 1_G^*\Omega^1_{G/V}$ et en dualisant, cela donne que  
$Lie(G)\simeq 1_G^*\TT_{G}$, de sorte que, si $G$ est lisse,
 $Lie(G)$ est un $V$-module libre de rang fini, dont la formation commute aux changements
de base. Le $K$-espace vectoriel $Lie(G)_K=Lie(G)\ot_V K$ est une algèbre de Lie sur $K$.

On introduit aussi $$Dist(G)_n=Hom_V(V[G]/I_1^{n+1},V)$$ et
$$Dist(G)= \varinjlim_n Dist(G)_n,$$ qui est une algèbre (cf I 7. de 
\cite{Jantzen}). On peut montrer, mais cela ne nous sera pas utile ici, que 
cette algèbre coïncide avec la fibre en $1_G$ du faisceau des opérateurs 
différentiels $\DD^{(0)}_G$. Après extension de $V$ à $K$, cette algèbre est 
l'algèbre enveloppante de $Lie(G)_K$. Le module $Lie(G)$ est canoniquement un sous-module 
de $Dist(G)$ et cet homomorphisme est un homomorphisme d'algèbres de Lie après extension
des scalaires à $K$.  

Soient $X(T)$ le groupe des caractères de $T$ ($X(T)=Hom(T,{\bf G}_m)$) et $Y(T)$ le groupe des 
sous-groupes de rang $1$ de $T$ ($Y(T)=Hom({\bf G}_m,T)$). Ce sont deux $\Ze$-modules libres de
rang fini et on dispose du crochet de dualité $<,>$: $X(T)\times Y(T)\rig \Ze$. 
En effet, soient $(\lambda,\mu)\in X(T)\times Y(T)$, alors 
$\lam\circ \mu $ définit un élément de $Hom_{Gr}(\Gm,\Gm)\simeq \Ze$. 
A un élément $\lam$ de $X(T)$, on associe un élément inversible de $V[T]$,
que l'on notera aussi $\lam$. 
Si $M$ est un $T$-module et $\lam\in X(T)$, et si $\Delta_m$ est l'application de
co-module $M\rig M\ot_V V[T]$, on note  
$M_{\lam}=\left\{m\in M\, |\, \Delta_M(m)=m\ot \lam\right\}.$
L'action par conjugaison de $T$ sur $Lie(G)$  notée $Ad$ se décompose comme d'habitude
$$Lie(G)=Lie(T)\bigoplus_{\alpha\in R}Lie(G)_{\alpha }.$$
Par définition, $R\subset X(T)$ est l'ensemble des racines de $G$. A chaque $\al$, on associe un 
élément $\al^{\vee}$ de $Y(T)$. L'ensemble des éléments $\alpha^{\vee}$ est noté $R^{\vee}$. Le quadruplet 
$(X(T),R,Y(T),R^{\vee})$ associé à la bijection $R\rig R^{\vee}$ et à l'accouplement 
$<,>$, constitue la donnée de racines de $G$. Comme $G$ est semi-simple, le couple 
$(R,X(T)\ot_{\Ze} K)$ est un système de racines sur le corps $K$ (cf chap. 6 de
\cite{bourbaki_lie_456}).

\subsubsection{Système de racines de l'algèbre de Lie d'un groupe algébrique}
\label{subsubsection-comp_syst_rac}
Soit $\lam\in X(T)$. L'application tangente $d\lam$ définit une application
$Lie(T)\rig Lie(\Gm)$ et donc un élément de $Lie(T)^*$, une fois fixée une coordonnée 
$t$ de $\Gm$. De même l'application tangente d'un élément $\mu$ de 
$Y(T)$ est une application $d\mu$ : $Lie(\Gm)\rig Lie(T)$. Après avoir 
identifié $Hom_{Lie}(Lie(\Gm),Lie(T))$ à $Lie(T)$, on peut voir l'élément $d\mu$ comme
un élément de $Lie(T)$, ce que nous ferons dans la suite. On obtient un 
accouplement canonique $\la d\lam,d\mu\ra$ en composant $d\lam\circ d\mu\in
Hom_{Lie}(Lie(\Gm),Lie(\Gm))\simeq V$. Par construction, on a, après choix d'une
coordonnée $t$ sur $\Gm$ : $\la\lam,\mu\ra=\la d\lam,d\mu\ra$ 
(car l'application tangente en $t=1$ de $t\mapsto t^n$ est la multiplication par $n$). 
L'application $X(T)\rig Lie(T)^*$ n'est pas injective en caractéristique $p>0$ 
(car le caractère $t\mapsto t^p$ est envoyé sur $0$), mais après tensorisation 
par $K$ on a un isomorphisme 
$X(T)\ot_{\Ze} K\simeq Lie(T_K)^*$ (resp. $X(T)\ot_{\Ze} k\simeq Lie(T_k)^*$).

Soit $\alpha\in R$, on lui associe $d\alpha\in Lie(T)^*$, qu'on notera $\alpha_*$. A la
co-racine $\alpha^{\vee}$, on associe de même $d\alpha^{\vee}\in Lie(T)$, qu'on notera 
$H_{\alpha}\in Lie(T)$. On note enfin $R_*=\left\{\alpha_*\,|\,\alpha\in R \right\}\subset 
Lie(T)^*$ et $R_*^{\vee}=\left\{H_{\alpha}|\,\alpha\in R \right\}\subset Lie(T)$. 

On définit comme d'habitude $s_{\alpha}$ : $X(T)\rig X(T)$ par 
$s_{\alpha}(\lam)=\lam - \la\lam,\alpha^{\vee}\ra\alpha$. Les applications 
$s_{\alpha}$ sont les réflexions associées au système de racines sur $K$
$(R,X(T)\ot_{\Ze}K)$. On définit de façon analogue pour $h\in Lie(T)^*$, $\alpha\in
R$, $\sigma_{\alpha_*}(h)=h - \la h, \alpha_*^{\vee}\ra\alpha_*$. On vérifie facilement
que les applications $\sigma_{\alpha_*}$ sont des réflexions. De plus, on a l'égalité,
pour $\alpha,\beta\in R$,
$\sigma_{\alpha_*}(\beta_*)=(s_{\alpha}(\beta))_*$ de sorte que les réflexions $\sigma_{\alpha_*}$
préservent $R_*$. On identifie ainsi les systèmes de racines $(R,X(T)\ot_{\Ze} K)$ et 
$(R_*,Lie(T_K)^*)$. 

On remarquera enfin que si $\lam\in X(T)$ et si $H\in Lie(T)$, alors $d\lam (H)=H(\lam)$. 
C'est vrai même si le tore $T$ n'est pas déployé. Pour voir cela, on commence par se ramener 
au cas où le tore est déployé après une extension fidèlement plate de la base. 
On est alors ramené à montrer cette égalité pour $\lam$ un caractère de $\Gm$, or, dans ce
cas, si $\lam(t)=t^n$, et $H=\der_1$ définie par $\der_1(f)=(\der f/\der t)(1)$, 
on a $d\lam (\der_1)=\der_1 (\lam)=n$.

Il nous reste maintenant à vérifier que le système de racines 
$(R_*,Lie(T_K)^*)$ est le système de racines associé à $Lie(G)$, c'est-à-dire que les
racines ainsi obtenues sont celles données par la représentation adjointe de $Lie(T)$ sur 
$Lie(G)$. Cela fait l'objet de la sous-section suivante. 

\subsubsection{Comparaison des actions adjointes}
L'action adjointe $Ad$ de $T $ sur $Lie(G)$ induit une action adjointe $ad$ de $Lie(T)$ sur $Lie(G)$ (I 7.11 de 
\cite{Jantzen}), et, si $\al\in R$, l'action de $H\in Lie(T)$ sur $v\in Lie(G)_{\al}$
est donnée par $H(v)=H(\al)v$. 
Dans la suite, nous vérifions que cette 
action de $Lie(T)$ déduite de l'action de conjugaison par $T$ 
est bien l'action adjointe induite par le crochet de Lie de $Lie(G)$. C'est
classique sur le corps des nombres complexes mais nous n'avons pas trouvé de réference 
sur une base plus générale. En particulier, cette assertion est valable sur un corps $k$ de 
caractéristique $p>0$ (dans ce cas, le lecteur identifiera $V=k=K$ dans les notations).
Pour ce faire, nous utilisons des formules sur les algèbres de distibutions citées dans
\cite{Jantzen}. Nous en déduirons la description de l'action de $Lie(T)$ sur $Lie(G)$ en
vertu des injections $Lie(T)\subset Dist(T)$, $Lie(G)\subset Dist(G)$ et $Dist(T)\subset
Dist(G)$. 

Si $H\in Dist(G)$, on note $\sigma'_G(H)=H\circ \sigma_G$. 
En I 7.18 (formule (3)) de \cite{Jantzen}, Jantzen
donne le calcul de l'action obtenue de $Dist(T)$ sur $Dist(G)$. L'immersion
diagonale $T\hrig T\times T$ définit sur $Dist(T)$ une co-multiplication 
$\Delta'_T$:$Dist(T)\rig Dist(T)\ot_V Dist(T)$ définie par $\Delta'_T(H)=1\ot H + H \ot
1$. Pour des raisons de fonctorialité, la co-multiplication sur $Dist(T)$ (resp. $\sigma'_T$) est la
restriction à $Dist(T)$ de la co-multiplication sur $Dist(G)$ (resp. $\sigma'_G$). Si $\Delta'_G(H)=\sum_i
H_i\ot H'_i$, Jantzen établit que $ad(H)H'=\sum_i H_iH'\sigma'_G(H'_i)\in
Dist(G)$. Comme $H\in Dist(T)$, 
on trouve $ad(H)(H')=HH'+H'\sigma'_T(H)$. Remarquons maintenant que 
$\sigma'_T(H)=-H$.  En effet, après extension fidèlement plate de la base, on se ramène au
cas où le tore $T$ est déployé, et finalement on est ramené au cas où
$T=\Gm=\spec(V[t,t^{-1}])$. Dans ce cas, $Dist(T)$ est libre de base $\der_1$ défini 
par $\der_1(f)=(\der f/\der t)(1)$.
  On calcule alors $$\sigma'_T(\der_1)(f)=\frac{\der \left(f(t^{-1})\right)}{\der t}(1)
=\left(-t^{-2}\left(\frac{\der f}{\der t}\right)(t^{-1})\right)(1)=-\der_1(f).$$ 
Finalement, cela montre que $ad(H)(H')=[H,H']$, de sorte que la décomposition 
de $Lie(G)=Lie(T)\bigoplus_{h_{\al}\in R}Lie(G)_{\al}$ est la décomposition définissant le 
système de racines de $Lie(G_K)$.  

\vspace{+3mm}
Un corollaire de toutes ces vérifications est qu'on peut identifier le 
 système de racines $(R, X(T_K))$ utilisé par Jantzen dans \cite{Jantzen}
et le système de racines $(R_*,Lie(T_K))$. Il reste 
à choisir un système de racines positifs. On remarquera que Kashiwara dans ~\cite{Kashi-D_flag}
et Beilinson-Bernstein dans ~\cite{BeBe} prennent la convention opposée. Nous adopterons ici la convention de 
Beilinson-Bernstein (et de Jantzen) en imposant que les racines venant du groupe de Borel constituent un système de 
racines négatives. 

\subsubsection{Poids dominants et réguliers}
Soit $\lam\in X(T)$, on note $\LL(\lam)$ le faisceau inversible associé à $\lam$ (I 5. de
\cite{Jantzen}). Le caractère $\lam$ induit un poids toujours noté $\lam$ sur $G_K$ et sur
$G_k$. Soient $V_{\lam}$, resp. $K_{\lam}$ et resp. $k_{\lam}$ les représentations de $T$, 
resp. $T_K$ et $T_k$ associées à $\lam$. 
A ces représentations, on associe les faisceaux $\LL(\lam)$ sur $X$, resp. $\LL(\lam_k)$ sur $X_k$ et $\LL(\lam_K)$ sur 
$X_K$. Il résulte de I 5.17 de \cite{Jantzen} que $i^*\LL(\lam)\simeq \LL(\lam_k))$
(resp. $j^*\LL(\lam)\simeq \LL(\lam_K)$). D'autre part, le faisceau $\LL(\lam)$ est
inversible sur $X$ d'après I 5.16 de \cite{Jantzen}.

On rappelle les définitions suivantes (II 2.6 de \cite{Jantzen}).
\begin{sousdef} Un poids $\lam $ de $G_K$ est appelé dominant (resp. régulier) si $\forall \alpha^{\vee}\in
R^{\vee}$, on a $\la\lam,\alpha^{\vee}\ra\geq 0$ (resp. $\forall \alpha^{\vee}\in
R^{\vee}$, on a $\la\lam,\alpha^{\vee}\ra >0$). 
\end{sousdef} 
De façon équivalente, $\lam $ est dominant si et seulement si le module de sections
globales $H^0(X,\LL_K(\lam))\neq 0$ et $\lam $ est régulier si et seulement si le faisceau 
$\LL_K(\lam)$ est ample. Si $\lam $ est dominant, $H^0(X,\LL_K(\lam))$ est la
représentation irréductible de $G_K$ de plus haut poids $\lam$ (à isomorphisme près).
Nous aurons besoin d'une variante sur $X$ du théorème de Kempf (II 4.5 de \cite{Jantzen}). 
\begin{sousprop}\label{paragraph-acyc_dom} 
Si $\lam$ est dominant, alors $\forall n\geq 1$, $H^n(X,\LL(\lam))=0$.
\end{sousprop} 
Le résultat est classique pour $\LL(\lam_K)$ et $\LL(\lam_k)$ sur $X_K$ et $X_k$
respectivement. Il résulte du lemme précédent que $i^*\LL(\lam)$ et $j^*\LL(\lam)$ sont
acycliques pour le foncteur sections globales. Comme la cohomologie commute à la limite
inductive sur un schéma noetherien, pour tout $n\in\Ne$, $H^n(X_K,\LL(\lam_K))=K\ot_V 
H^n(X,\LL(\lam))$. Comme ces groupes sont nuls pour $i\geq 1$, on voit que les groupes 
$H^n(X,\LL(\lam))$ sont de torsion et donc de torsion finie car ce sont des 
$V$-modules de type fini par les théorème généraux. La longue suite exacte de cohomologie
associée à la suite exacte courte 
$$0\rig \LL(\lam)\sta{\cdot \pi}{\rig}\LL(\lam)\rig i_*\LL(\lam_k)\rig 0,$$
donne des surjections, pour tout $n\geq 1$ 
$$\cdot\pi\,\colon\,H^n(X,\LL(\lam))\trig H^n(X,\LL(\lam)),$$
ce qui montre finalement que ces groupes sont nuls puisqu'ils sont de torsion finie.
\subsection{Faisceau tangent sur un espace homogène}
On termine par le fait classique suivant.
\begin{prop}\label{subsubsection-TX}
Le faisceau $\TT_X$ est engendré par ses sections globales.
\end{prop} 
L'action à gauche de $G$ sur $X$ (resp. $G$) munit $\TT_X$
 (resp. $\TT_G$) d'une structure de $G$-module équivariant (chapitre 1 de \cite{GIT}). En particulier, 
le groupe abélien $\Ga(G,\TT_G)$ est un $G$-module. 
Le module des dérivations
invariantes $\Ga(G,\TT_G)^G$ s'identifie à $Lie(G)$ (II 4 6.5 de \cite{Demazure_gr_alg}).
La formation de ce module commute donc aux changements de base. D'après 
II 4 6.3 de \cite{Demazure_gr_alg}, on dispose d'une application canonique 
$Lie(G)\rig \Ga(X,\TT_X)$. On en déduit une application 
canonique $u$ : $\OO_X\ot_{V}\Ga(G,\TT_G)^G \rig \TT_X$.

Dans le cas où la base est un corps algébriquement clos, il est bien connu que 
$u$ est surjectif et cela résulte 
du fait que l'action de $G$ sur $X$ est transitive. Rappelons la 
démonstration dans ce cas. Supposons maintenant que 
$S=\spec \overline{l}$ où $\overline{l}$ est un corps algébriquement clos.  Notons  
 $e=1_GP$, la fibre $i_e^*\TT_X$ est un quotient de $Lie(G)$ d'après 
II 4.2 de \cite{Jantzen}. Il suffit de montrer que $u$ est surjectif au-dessus des points 
fermés d'après le lemme de Nakayama. Soit $x\in X(\overline{l})$, alors il existe $g\in G(\overline{l})$ 
tel que $x=gP$. On dispose alors d'opérateurs de translation $\lam_g$: 
$k(x)=\overline{l}\simeq k(e)=\overline{l}$ et 
$\rho_g$ : $ i_x^*\TT_X\simeq  i_e^*\TT_X,$ semi-linéaires par rapport aux 
$\lam_g$ tels que le diagramme suivant soit commutatif
$$\xymatrix{\Ga(G,\TT_G)^G\ar@{->}[r]^(.6){u(x)} \ar@{=}[d] & 
i_x^*\TT_X \ar@{->}[d]^{\rho_g}_{\wr}\\
        \Ga(G,\TT_G)^G \ar@{->>}[r]^(.6){u(e)} & i_e^*\TT_X .}$$
Cela montre que $u(x) $ est surjectif et donc finalement que $u$ est surjectif. 
Sur une base $S$ générale, il suffit de montrer la surjectivité de $u$ en tout point fermé 
$s$ de $S$, d'après le lemme de Nakayama. Soient $i_s$ l'immersion fermée correspondante 
à $s$, $k(s)$ le corps résiduel de $s$ et 
$\overline{l}$ une clôture algébrique de $k(s)$. Après 
application de $i_s^*$, l'application $u$ donne une application
$u_s$: $\OO_{X_s}\ot_{k(s)}\Ga(G_s,\TT_{G_s})^{G_s}\rig \TT_{X_s}$. 
Après extension des scalaires à $ \overline{l}$, cette flèche est surjective d'après ce 
qui précède. Par fidèle platitude de $\overline{l}$ sur $k(s)$, la flèche $u_s$ est surjective
et donc $u$ est surjectif.

\vspace{+3mm} Dans la sous-section suivante, on explique pourquoi il suffit de montrer le 
théorème de Beilinson-Bernstein après extension fidèlement plate de la base.
\subsection{Changements de base fidèlement plats}
Dans cette sous-section, $X$ est un schéma lisse sur $S$, dont le complété formel est 
$\XX$. 
Soit $\DD_{\XX}$ l'un des faisceaux $\Dcm_{\XX,\Qr}$ pour un certain $m$ ou $\Ddag_{\XX,\Qr}$.
Soient $V'$ un anneau de valuation discrète d'inégales caractéristiques $0,p$, qui 
est une $V$-algèbre finie fidèlement plate, $S'=\spec V'$, $X'=X\times_S S'$,
$\SS'=\spf V'$, $\XX'=\SS'\times_{\SS}\XX$. On dispose alors de la proposition suivante
\begin{prop} \label{subsubsection-ext_fplat}
Si $\XX'$ est $\DD_{\XX'}$-affine, alors $\XX$ est
$\DD_{\XX}$-affine (resp. si $X'$ est $\DD_{X'}$-affine, alors $X$ est $\DD_X$-affine).
\end{prop} 
\begin{proof} Soit $\MM$ un $\DD_{\XX}$-module cohérent. Il suffit de montrer qu'il est
acyclique pour le foncteur $\Gamma$ et engendré par ses sections globales comme 
$\DD_{\XX}$-module. Le $\DD_{\XX'}$-module $V'\ot_V\MM$ est cohérent, donc acyclique pour
le foncteur $\Ga$ et engendré par ses sections globales. De plus, par platitude du morphisme $V\rig V'$, 
on a pour tout $n\geq 0$, 
$$H^n(\XX',V'\ot_V\MM)=V'\ot_V H^n(\XX,\MM),$$
de sorte que $\forall n\geq1$, les groupes $H^n(\XX,\MM)$ sont nuls.
De façon analogue, on a une surjection 
$$\Ddag_{\XX'}\ot_{\Ga(\XX',\Ddag_{\XX'})}\Ga(\XX,V'\ot_V \MM)\trig V'\ot_V \MM,$$ 
dont on déduit que $\MM$ est engendré par ses sections globales par fidèle platitude de 
$V'$ sur $V$. 
\end{proof}
{\bf Remarque}. On donnera une réciproque à cet énoncé en ~\ref{subsubsection-ext_plat} en supposant seulement
que $V\rig V'$ est fini et plat, pour les espaces homogènes, 
et plus généralement les schémas vérifiant l'hypothèse (H) de ~\ref{section-hypH}.


\vspace{+3mm}
Passons maintenant à la démonstration du théorème principal.

\section{Un critère pour passer du cas algébrique au cas formel}
\label{section-hypH}
Le théorème d'annulation va provenir d'un énoncé plus général sur des schémas projectifs 
sur $S$, sur lesquels le faisceau structural est acyclique et dont le faisceau tangent est 
engendré par ses sections globales. Il s'agit de donner un critère pour passer d'un
énoncé d'acyclicité pour les $\DD$-modules sur un schéma projectif lisse $X_K$ sur un corps 
$p$-adique à un tel énoncé d'acyclicité pour les $\Ddag_{\XX,\Qr}$-modules cohérents sur le complété 
formel d'un modèle entier de $X_K$.
Dans cette partie, on suppose que $X$ est un schéma projectif lisse vérifiant les
hypothèses suivantes (H):\be \item[(i)]Le faisceau $\OO_X$ est acyclique pour le
foncteur $\Ga$.
\item[(ii)]Le faisceau tangent $\TT_X$ est engendré par ses sections globales.\ee 
Soit $\LL$ un faisceau inversible sur $X$. On note encore $\LL$ le faisceau inversible 
obtenu à partir de $\LL$ sur la complétion formelle $\XX$ de $X$ et pour n'importe quel 
symbole $\sharp$ égal à $(m)$ ou $\dagger$, on introduit comme en ~\ref{subsection-diff_finv}
le faisceau des opérateurs différentiels $\DD^{\sharp}_{\XX}(\LL)$ à valeurs dans $\LL$
défini par $\DD^{\sharp}_{\XX}(\LL)=\LL \ot_{\OO_X}\DD^{\sharp}_{\XX}
\ot_{\OO_X}\LL^{-1}$. 

L'un des points clefs de la démonstration consiste à avoir un résultat
de finitude de la torsion des groupes $H^n(X,\Dm_{X}(s))$ pour $s\in\Ze$ fixé et 
$n\geq 1$. Pour
cela on utilise les techniques de $\cite{Hu1}$ et les techniques de Kashiwara exposées 
en 1.4 de \cite{Kashi-D_flag}. L'idée consiste à donner un critère analogue à 
celui de Kashiwara, à torsion finie près.

\vspace{+3mm}
Fixons maintenant un faisceau ample inversible $\OO_X(1)$ sur $X$. Tout 
$\OO_X$-module cohérent est quotient d'un faisceau du type 
$\OO_X(-r)^a$, avec $a,r\in\Ne$. De plus, il existe $U\in\Ne$, tel que 
 pour tout $ u\geq U$, le faisceau $\OO_X(u)$ est engendré par ses sections globales au 
sens où la flèche suivante est surjective $$ \OO_X\ot_V\Ga(X,\OO_X(u))\trig \OO_X(u),$$
et est acyclique pour le foncteur $\Ga$.

Ces propriétés ainsi que les propriétés cohomologiques 
du faisceau structural $\OO_X$ seront essentielles dans ce qui suit.
Si $\EE$ est un $\OO_X$-module, $\EE(s)$ pour $s\in\Ze$ désigne 
$$\EE(s)=\EE\ot_{\OO_X}(\OO_X(1))^{\ot s}.$$

Soit $\LL$ un $\OO_X$-module localement libre de rang $1$
(induisant $\LL$ sur $X_K$ et $\LL$ sur $\XX$).
L'objet de cette section est de montrer le théorème suivant.
\begin{surthm} \label{subsection-thm-int}Soit $X$ un schéma lisse vérifiant les hypothèses (H), $\XX$ le schéma 
formel associé, on a l'énoncé suivant:
si $X_K$ est $\DD_{X_K}$-affine (resp. $\DD_{X_K}(\LL)$-affine), alors $\XX$ est
$\DD^{\dagger}_{\XX,\Qr}$-affine
(resp. $\DD^{\dagger}_{\XX,\Qr}(\LL)$-affine).
\end{surthm} 
Ce théorème sera démontré en ~\ref{subsubsection-dem_thm_int} et en ~\ref{subsubsection-eq_cat}. 
On montrera aussi que les hypothèses entraînent que $\XX$ est $\Dcm_{\XX,\Qr}$-affine 
pour tout entier $m$. Pour un $m$ fixé, la démonstration repose sur la
structure de l'algèbre graduée $\gr\pg\Dm_X$ et sur les résultats classiques de la
cohomologie des faisceaux cohérents sur un schéma projectif. On remarquera que 
si $\LL$ est un $\OO_X$-module inversible, 
$$\gr\pg\Dm_X(\LL)\simeq \LL\ot_{\OO_X}\gr\pg\Dm_X\ot_{\OO_X}\LL^{-1},$$
et donc $\gr\pg\Dm_X(\LL)\simeq \gr\pg\Dm_X$ puisque l'algèbre graduée 
$\gr\pg\Dm_X$ est commutative. Grâce à cette remarque, le lecteur se rendra compte 
que la démonstration du théorème est la même dans le cas de $\Dcm_{\XX,\Qr}$ (resp.
$\Ddag_{\XX,\Qr}$) et de 
$\Dcm_{\XX,\Qr}(\LL)$ (resp. $\Ddag_{\XX,\Qr}(\LL)$). 
Pour éviter d'alourdir les notations, nous ferons la démonstration pour le cas de 
 $\Dcm_{\XX,\Qr}$ (resp. $\Ddag_{\XX,\Qr}$).
%
\subsection{Résultats à un niveau fini}
%
Dans cette partie, $m$ est fixé.
Le premier résultat consiste à établir que si $\MM$ est un $\Dm_X$-module cohérent, alors 
$\MM(r)$ est acyclique pour le foncteur $\Ga(X,.)$ pourvu que $r$ soit assez grand. Ce
résultat repose sur un résultat analogue pour les modules sur l'algèbre graduée 
$\gr\pg\Dm_X$ et pour les $\OO_X$-modules. Pour l'algèbre graduée $\gr\pg\Dm_X$, nous avons en 
effet la proposition suivante.
\begin{prop} \label{subsubsection-grD}
Il existe $r_0\in\Ne$ tel que $\forall r\geq r_0,\forall n\geq 1$,
$H^n(X,\gr\pg\Dm_X(r))=0$.
\end{prop} 
En particulier, pour tous $r\geq r_0$, $n\geq 1$, $t\geq 0$, les groupes
$H^n(X,\gr_t\Dm_X(r))$ sont nuls.
\begin{proof} D'après (H), il existe une surjection $\OO_X^a\rig\TT_X$, d'où on déduit 
(\ref{subsubsection-devis}) un morphisme surjectif
de faisceaux cohérents d'algèbres graduées $$\CC=\bSm(\OO_X^a)\trig \bSm(\TT_X),$$ et 
cette dernière algèbre graduée s'identifie à l'algèbre graduée $gr\pg\Dm_X $
(\ref{subsubsection-isom_alg_grad}). Il suffit donc de montrer l'assertion pour un 
$\CC$-module cohérent $\EE$. Comme 
$X$ est noetherien, et $\EE$ est un $\OO_X$-module quasi-cohérent, 
$\EE$ est limite inductive de ses sous $\OO_X$-modules cohérents
$\EE_i$ pour $i\in I$.
Comme $\EE$ est un $\CC$-module cohérent et que le faisceau d'algèbres $\CC$ est 
à sections noetheriennes sur les ouverts affines, il existe une surjection $\CC$-linéaire 
$$\CC\ot_{\OO_X}\EE_i\trig \EE,$$ et donc une surjection $\CC$-linéaire 
$$\CC(s_0)^{a_0}\trig \EE,$$ avec $a_0\in\Ne$ et $s_0\in\Ze$. Si $\EE$ est gradué, on
peut faire en sorte que la surjection soit graduée mais cela ne sera pas important ici. 
Finalement, on peut construire ainsi de proche en proche, une résolution de longueur 
$N+1$ de $\EE$ par des $\CC$-modules cohérents, du type suivant
$$0\rig \FF \rig \EE_{N-1} \rig \ldots \EE_1\rig\EE_0 \rig 0$$ avec 
$$\EE_i=\CC(s_i)^{a_i},$$ pour $0\leq i \leq N-1$, $a_i\in\Ne$, $s_i\in\Ze$. Tensorisons cette résolution par 
$\OO_X(r)$ avec $r\geq r_0=max\{U,U-s_i\}_{0\leq i \leq N-1}$. Les termes d'indice 
$0,\ldots,N-1$ de cette 
résolution sont des modules du type $\CC(r+s_i)^{a_i}$. Ces modules sont sommes directes 
de composantes homogènes du type $\OO_X(r+s_i)^{b_i}$. D'après notre choix de $r_0$ et de
$U$, ces modules sont acycliques pour le foncteur $\Ga(X,.)$.  
On utilise alors l'énoncé suivant: 
si un $\CC$-module cohérent $\EE$ admet une résolution de longueur $\geq N+1$, 
$$ 0\rig \FF \rig \EE_{N-1} \rig \ldots \EE_1\rig\EE_0 \rig 0,$$
par des $\CC$-modules cohérents, telle que pour un entier 
 $0\leq i\leq N-1$, les modules $\EE_0,\EE_1,\ldots ,\EE_i$
sont acycliques pour le foncteur $\Ga(X,.)$, 
alors, pour tout $N-i \leq n \leq N$, $H^n(X,\EE)=0$. En effet, 
on dispose d'une suite spectrale bi-régulière associée à cette résolution
$$ H^j(X,\EE_t)\Longrightarrow H^{j-t}(X,\EE),$$
et, par hypothèse, pour $n\geq N-i$, les seuls termes $H^j(X,\EE_t)$ 
intervenant dans le gradué de la
filtration sur $H^{n}(X,\EE)$ vont correspondre à des valeurs de $t\leq i$ et $j\geq 1$,
pour lesquelles les groupes $H^{j-t}(X,\EE)$ sont nuls. Si bien que l'aboutissement $H^n(X,\EE)$ 
est égal à $0$ pour $N-i \leq n \leq N$. \end{proof}

On en déduit le corollaire 
\begin{cor}\label{subsubsection-Dcoh} \be \item[(i)] $\forall r\geq r_0, \forall n\geq 1,
\forall t\geq 0, \, H^n(X,\Dm_{X,t}(r))=0 ,$
\item[(ii)] $ \forall r\geq r_0, \forall n\geq 1,
 \, H^n(X,\Dm_{X}(r))=0 ,$
\item[(iii)] Pour tout $\Dm_X$-module cohérent $\MM$, il existe 
$r_1\in\Ne$, tel que $\forall r\geq r_1$, $\forall n\geq 1$, $H^n(X,\MM(r))=0.$
\ee
\end{cor} 
\begin{proof} Le (ii) résulte du (i) par passage à la limite inductive. On montre le (i) par
récurrence sur $t$. Pour $t=0$, $\Dm_{X,0}(r)=\OO_X(r)$, qui est acyclique pour le foncteur 
$\Ga$ car $r\geq r_0\geq U$. Pour tout $t\geq 1$, et
tout $r\geq r_0$, on dispose de suites exactes courtes 
$$0\rig \Dm_{X,t-1}(r) \rig \Dm_{X,t}(r) \rig \gr_t\Dm_{X}(r)\rig 0.$$
Comme le faisceau $\gr_t\Dm_{X}(r)$ est acyclique pour $\Ga$ d'après la proposition 
précédente, on voit par récurrence sur $t$ qu'il en est de même pour les faisceaux 
$\Dm_{X,t}(r)$ après application de la suite exacte longue de cohomologie pour 
$\Ga$. 

Pour le (iii), on remarque, en procédant comme en ~\ref{subsubsection-grD}, que $\MM$ admet une résolution
du type suivant
$$0\rig \EE_N \rig \EE_{N-1} \rig \ldots \EE_1\rig\EE_0 \rig 0$$ avec 
$$\EE_i=\Dm_X(s_i)^{a_i},$$ pour $0\leq i \leq N-1$ et $s_i\in\Ze$. Tensorisons cette 
résolution par $\OO_X(r_1)$ avec $r_1=max \{r_0+s_i\}_{0\leq i \leq N-1}$. On voit que le (iii) résulte du (ii) en
utilisant le même argument de suite spectrale qu'en ~\ref{subsubsection-grD}.\end{proof}

\vspace{+3mm}
Dans la suite de cette sous-section, on se place sous les hypothèses du théorème 
et on suppose que $X_K$ est $\DD_{X_K}$-affine.

Le point clé de la démonstration est le résultat suivant, qui concerne la cohomologie des 
faisceaux $\Dm_X(s)$ pour $s\in\Ze$.
\begin{prop}\label{subsubsection-propcle}
Soit $s\in\Ze$. Alors $\forall n\geq 1$, $H^n(X,\Dm_X(s))$ est un groupe
de torsion finie. 
\end{prop} 
\begin{proof}
Le faisceau $\Dm_{X,\Qr}(s)$ est un $\DD_{X_K}$-module cohérent, donc pour $n\geq 1$ les 
groupes $H^n(X_K,\Dm_{X,\Qr}(s))$ sont nuls par hypothèse, et sont égaux à 
$H^n(X,\Dm_X(s))\ot_V K$, par commutation de la cohomologie à la limite inductive, 
de sorte que les groupes $H^n(X,\Dm_X(s))$ sont de torsion pour $n\geq 1$. Pour voir que
la torsion est finie, on s'inspire des arguments de 1.4 \cite{Kashi-D_flag}.

Fixons $u\geq max\{r_0-s,U\}$. Commençons par remarquer qu'on a une section $\ti{\tau}$, 
$\DD_{X,\Qr}$-linéaire à gauche, à la surjection canonique $\ti{\sigma}$
$$\xymatrix {\DD_{X,\Qr}\ot_{\OO_{X_K}}\OO_{X_K}(-u)\ot_K\Ga(X,\OO_X(u))\ar@{->>}[r]^(.8){\ti{\sigma}} 
& \DD_{X,\Qr}\ar@{-->}@<5pt>[l] ^(0.2){\ti{\tau}} .}$$
On part en effet de la surjection canonique 
$$\DD_{X,\Qr}\ot_K \Ga(X_K,\OO_{X_K}(u))\trig \DD_{X,\Qr}\ot_{\OO_{X_K}}\OO_{X_K}(u).$$
Comme le foncteur $\Ga$ est exact sur la catégorie des $\DD_{X,\Qr}$-modules cohérents, on trouve une
surjection
$$\Ga(X,\DD_{X,\Qr}(-u))\ot_K \Ga(X,\OO_{X_K}(u))\trig \Ga(X,\DD_{X,\Qr}),$$
dont on trouve une section $\ti{\tau}$, $\DD_{X,\Qr}$-linéaire à gauche, en relevant 
$1\in\Ga(X,\DD_{X,\Qr})$. Observons maintenant qu'il existe $i\in\Ne$ tel que 
$\pi^i\ti{\tau}(1)\in \Ga(X,\Dm_{X}(-u))\ot_V \Ga(X,\OO_{X}(u))$. Définissons 
$\overline{\tau}$ comme l'unique application $\Dm_X$-linéaire 
$\Dm_X \rig \Dm_X\ot_{\OO_X}\OO_X(-u)\ot_V\Ga(X,\OO_X(u))$ définie par 
$\overline{\tau}(1)=\pi^i\ti{\tau}(1)$. 
En particulier, le conoyau de l'application suivante $\ti{\sigma}_m$ est annulé par $\pi^i$ 
$$\Dm_X\ot_V \Ga(X,\OO_{X}(u))\rig \Dm_X(u).$$
Par construction, on a alors 
$$\ti{\sigma}_m\circ \overline{\tau}_m=\pi^i id_{\Dm_X}.$$
En dualisant (après application de $\HH om_{\Dm_X}(.,\Dm_X)$) et en tensorisant 
par $\OO_X(s)$, on trouve le 
diagramme suivant de $\Dm_X$-modules à droite 
$$\xymatrix{ \Dm_X(s) \ar @{^{(}->} [r] ^(.3){\sigma_m} & \Dm_X(u+s)\ot_V \Ga(X,\OO_{X}(u))^* 
\ar@{-->}@<5pt>[l] ^(0.7){\tau_m},}$$
et on a la relation $$\tau_m\circ \sigma_m=\pi^i id_{\Dm_X(s)}.$$
Considérons maintenant le diagramme commutatif suivant, pour $n\geq 1$, $t\geq 0$,
$$\xymatrix{ H^n(X,\Dm_{X,t}(s))\ar @{->}[r]\ar @{->}[d]_{a_t} \ar @{} [dr] |{\circlearrowleft}
& H^n(X,\Dm_{X,t}(u+s)\ot_V\Ga(X,\OO_X(u))^*)\ar @{->}[d]\\
             H^n(X,\Dm_{X}(s))\ar @{->}[r]_(.3){b} &
H^n(X,\Dm_{X}(u+s)\ot_V\Ga(X,\OO_X(u))^*),}$$
avec $b=H^n\circ \sigma_m$. Posons aussi $c=H^n\circ \tau_m$, de sorte que 
$$c\circ b=\pi^i id_{H^n(X,\Dm_X(s))}.$$ 
On a choisi $u$ pour que les deux termes de la colonne de droite du diagramme soient nuls,
ce qui implique que $b\circ a_t=0$ et donc, en composant avec $c$, que $\pi^ia_t=0$. En
passant à la limite inductive sur $t$, cela nous donne que, pour $n\geq 1$ fixé, 
$\pi^iH^n(X,\Dm_{X}(s))=0$ et donc l'énoncé de la proposition.\end{proof}

On en tire le corollaire suivant.
\begin{cor} \label{subsubsection-res_fini}Soit $\MM$ un $\Dm_X$-module cohérent, alors $\forall n\geq 1$, 
$H^n(X,\MM)$ est de torsion finie.
\end{cor} 
\begin{proof} Puisque $\Dm_X$ est à sections noetheriennes sur les affines, on peut 
procéder comme en ~\ref{subsubsection-grD} et $\MM$ admet une résolution $\Dm_X$-linéaire de 
longueur $\geq N$ par des modules sommes directe de modules du type $\Dm_X(s)$. En
procédant comme en ~\ref{subsubsection-grD} et en utilisant le même argument de suite
spectrale, on voit que les groupes $H^n(X,\MM)$ sont de torsion finie pour $n\geq 1$.
\end{proof}

\vspace{+1mm}
Contrairement aux résultats de la partie 4 de \cite{Hu1}, nous ne pouvons pas 
donner d'énoncé de finitude des sections globales de $\Dm_X$. Cela vient 
du fait qu'on ne sait pas si $\Ga(X,\gr\pg\Dm_X)$ est finie sur 
$\gr\pg\Ga(X,\Dm_X)$ en général. Sous les hypothèses (H), on a le résultat de finitude 
suivant sur la cohomologie des $\bSm(\TT_X)$-modules.
 Reprenons les notations de 
~\ref{subsubsection-grD}.
\begin{prop} L'algèbre $\Ga(X,\bSm(\TT_X))$ est noetherienne. De plus, si 
$\EE$ est un $\bSm(\TT_X)$-module cohérent, pour tout $n\in\Ne$, 
le module $H^n(X,\EE)$ est de type fini sur $\Ga(X,\bSm(\TT_X))$.
\end{prop} 
\begin{proof}
Fixons un plongement projectif $i'$ : $X\hrig Y=\bP_{V}^N$ tel que 
$\OO_X(1)={i'}^*\OO_Y(1)$. Par hypothèse (H), le faisceau $\bSm(\TT_X)$ est un 
$\CC$-module cohérent. Posons $\CC'=\bS^{(m)}(\OO_Y^a)$. Alors le faisceau 
$i'_*\bSm(\TT_X)$ est un $\CC'$-module cohérent. Le morphisme $i'$ est fini, de sorte
qu'il suffit de montrer que $C'=\Ga(Y,\CC')$ est une algèbre noetherienne et que si $\FF$ est un 
$\CC'$-module cohérent, les modules $H^n(Y,\FF)$ sont de type fini sur 
$C'$. Remarquons que $C'=\bS^{(m)}(V^a)$ est une 
$V$-algèbre de type fini (\ref{subsection-sym_m}) et est donc noetherienne. 
Par le même argument de suite spectrale et d'existence de résolutions particulières que 
celui qui est utilisé en ~\ref{subsubsection-grD}, il suffit de montrer que 
les modules $H^n(Y,\CC'(s))$ sont de type fini sur $C'$ pour tout $s\in\Ze$ et tout 
$n\in\Ne$. Or, par commutation de la cohomologie à la limite inductive, on a 
$$H^n(Y,\CC'(s))\simeq C'\ot_V H^n(Y,\OO_Y(s)).$$
Comme les groupes $H^n(Y,\OO_Y(s))$ sont des $V$-modules de type fini, cela 
donne le fait que les $C'$-modules $H^n(Y,\FF)$ sont de type fini pour tout 
$\CC'$-module cohérent $\FF$. C'est en particulier le cas pour 
$\Ga(X,\bSm(\TT_X))$, qui est donc une algèbre noetherienne. 
\end{proof}

\vspace{+3mm}
Il s'agit désormais de passer au cas du schéma formel $\XX$. 
\subsection{Passage au schéma formel}
On suppose dans toute cette sous-section que $X$ est un $S$-schéma vérifiant l'hypothèse
(H). On commence par montrer que la catégorie des 
$\Dcm_{\XX}$-modules cohérents est engendrée par les modules du type 
$\Dcm_{\XX}(-r)$ pour $r\in\Ze$. Cela correspond à la proposition 3.5 de \cite{Hu1}. 
La démonstration est identique et suit de ~\ref{subsubsection-Dcoh}.

\begin{prop}\label{subsubsection-dcm_gener} Soit $\MM$ un $\Dcm_{\XX}$-module cohérent
(resp. un $\Ddag_{\XX,\Qr}$-module cohérent). 
\item[(i)]Il existe 
$r_2\in\Ne$, tel que $\forall r\geq r_2$, $\forall n\geq 1$, $H^n(X,\MM(r))=0.$
\item[(ii)]Il existe $(a,r)\in\Ne^2$ et une 
surjection $\Dcm_{\XX}$-linéaire $\left(\Dcm_{\XX}(-r)\right)^a\trig \MM$
(resp. $(a,b,r,s)\in\Ne^4$ et une résolution à $2$ termes $\Ddag_{\XX,\Qr}$-linéaire 
$\left(\Ddag_{\XX,\Qr}(-s)\right)^b\rig\left(\Ddag_{\XX,\Qr}(-r)\right)^a\rig \NN\rig 0$).
\end{prop} 
Comme corollaire, on en déduit une réciproque à~\ref{subsubsection-ext_fplat}. Reprenons
les notations de cet énoncé en supposant seulement que le morphisme 
$V\rig V'$ est fini et plat. Alors, on a 
\begin{cor} \label{subsubsection-ext_plat}
Si $\XX$ est $\DD_{\XX}$-affine, $\XX'$ est $\DD_{\XX'}$-affine.
\end{cor} 
\begin{proof}Soit $\MM$ un $\DD_{\XX'}$-module cohérent. Il est clair que 
$X'$ vérifie (H). Soit $\MM$ un $\DD_{\XX'}$-module cohérent. 
Ce module admet une résolution de longueur arbitrairement grande par des modules 
du type $\DD_{\XX'}(-r)^a$. Par le même argument de suite spectrale qu'en 
~\ref{subsubsection-grD}, il suffit de montrer que les modules $\DD_{\XX'}(-r)$ sont acycliques pour 
$\Ga$ pour vérifier qu'il en est de même pour $\MM$. Or, 
$\DD_{\XX'}(-r)=V'\ot_V\DD_{\XX}(-r)$, et comme $V\rig V'$ est plat, on a, pour tout 
$n\geq 0$, $$H^n(\XX',\DD_{\XX'}(-r))=V'\ot_V H^n(\XX,\DD_{\XX}(-r)),$$
d'où l'énoncé d'acyclicité. De plus, le faisceau $\DD_{\XX'}(-r)$ qui est obtenu par
changement de base à partir de $\DD_{\XX}(-r)$, est engendré par ses sections globales 
comme $\DD_{\XX'}$-module. Comme le foncteur $\Gamma$ est exact pour les
$\DD_{\XX'}$-modules cohérents, on dispose d'un diagramme commutatif
$$\xymatrix{ \DD_{\XX'}\ot_{\Ga(\XX',\DD_{\XX'})}\Ga(\XX',\DD_{\XX'}(-r)^a)\ar @{->>}[r]
\ar @{->>}[d]&
\DD_{\XX'}\ot_{\Ga(\XX',\DD_{\XX'})}\Ga(\XX',\MM)\ar @{->}[d]\\
\DD_{\XX'}(-r)^a\ar @{->>}[r] & \MM ,}$$
 qui montre que $\MM$ est engendré par ses sections globales comme $\DD_{\XX'}$-module.
\end{proof}
A partir de maintenant, dans tout le reste de cette sous-section, on suppose 
que $X_K$ est $\DD_{X_K}$-affine. Grâce à la proposition 
précédente~\ref{subsubsection-dcm_gener}, on est 
ramené à contrôler les groupes $H^n(X,\Dm_{X}(s))$, en vue de l'énoncé d'acyclicité. 
Dans la partie 3. de ~\cite{Hu1}, on utilise le fait que  
les groupes $H^n(X,\Dm_{X}(s))$ sont des $V$-modules de type fini de torsion pour $n\geq
1$ et $s\in\Ze$. La différence ici 
est que l'on n'a pas de propriétés de finitude sur $V$, mais on sait que ces groupes sont
de torsion finie d'après ~\ref{subsubsection-propcle}. Cependant, le lecteur pourra vérifier 
que dans la partie 3 de \cite{Hu1}, seule la finitude de la torsion des groupes $H^n(X,\MM)$
est utilisée. Cette propriété permet de vérifier des conditions de Mittag-Leffler pour 
les groupes de cohomologie  $H^n(X_i,\Dm_{X_i}(s))$ 
pour $i$ variable, et permettent des passages à la limite pour la cohomologie.
Le résultat suivant se démontre comme la proposition 3.2 de~\cite{Hu1} compte tenu de
~\ref{subsubsection-propcle}.
\begin{prop}\label{subsubsection-coh_comp} Soit $\MM$ un $\Dm_{\XX}$-module cohérent et 
$\what{\MM}=\varprojlim_i \MM/\pi^{i+1}\MM$. Alors
\be\item[(i)]$\forall n\in\Ne,\, H^n(\XX,\what{\MM})=\varprojlim_i
H^n(X_i,\MM/\pi^{i+1}\MM),$
\item[(ii)] $\forall n\geq 1,\, H^n(\XX,\what{\MM})=H^n(X,\MM).$
\ee
\end{prop} 
En particulier, pour $n\geq 1$, $H^n(\XX,\what{\MM})$ est un $V$-module de torsion finie.

A partir du (ii) de la proposition~\ref{subsubsection-dcm_gener}, on peut procéder comme pour 
~\ref{subsubsection-res_fini}, ce qui donne le résultat de finitude suivant
\begin{prop} Soit $\MM$ un $\Dcm_{\XX}$-module cohérent, alors pour 
tout $n\geq 1$, les groupes $H^n(\XX,\MM)$ sont de torsion finie.
\end{prop} 
Donnons les conséquences de ces résultats pour les $\Dcm_{\XX,\Qr}$-modules cohérents.
Soit $\NN$ un $\Dcm_{\XX,\Qr}$-module cohérent. Comme l'espace topologique 
associé à $\XX$ est noetherien, il existe d'après 3.4.5 de \cite{Be1} un $\Dcm_{\XX}$-module cohérent $\MM$ tel que 
$\NN=\MM\ot_{V}K.$ Soit maintenant $\NN$ est un $\Ddag_{\XX,\Qr}$-module cohérent, d'après
3.6.2 de \cite{Be1}, il existe
un $\Dcm_{\XX,\Qr}$-module cohérent $\NN_0$ tel que 
$$\NN\simeq \Ddag_{\XX,\Qr}\ot_{\Dcm_{\XX,\Qr}}\NN_0.$$
Comme la cohomologie commute à la limite inductive sur $\XX$, les deux propositions 
précédentes nous 
permettent de montrer les énoncés suivants, en procédant comme en 3.5 de \cite{Hu1} et
montrent la partie \og acyclicité \fg ~de ~\ref{subsection-thm-int}.
\begin{prop}\label{subsubsection-dem_thm_int}
Soit $\NN$ un $\Dcm_{\XX,\Qr}$-module cohérent (resp. un
$\Ddag_{\XX,\Qr}$-module cohérent), alors 
 $\forall n\geq 1$, $H^n(X,\NN)=0,$
\end{prop} 
Indiquons comment passer à des énoncés de $\Dcm_{\XX,\Qr}$-affinité (resp. 
$\Ddag_{\XX,\Qr}$-affinité). Fixons $u\geq max\{r_0,U\}$ et reprenons
les applications $\tau_m$ et $\sigma_m$ pour $s=0$ construites lors de la démonstration de 
~\ref{subsubsection-propcle}. Complétons ces applications, tensorisons par 
$\OO_{\XX}(-u)$ et inversons $\pi$, cela nous donne
un diagramme 
$$\xymatrix{ \Dcm_{\XX,\Qr}(-u) \ar @{^{(}->} [r] ^(.4){\what{\sigma}_m} & \Dcm_{\XX,\Qr}\ot_V \Ga(X,\OO_{X}(u))^* 
  \ar@{-->}@<5pt>[l] ^(0.6){\what{\tau}_m},}$$
 et l'application $s_m=\pi^{-i}\what{\tau}_m$ est une section $\Dcm_{\XX,\Qr}$-linéaire 
de $\what{\sigma}_m$. On obtient une surjection $\Dcm_{\XX,\Qr}$-linéaire 
$\Dcm_{\XX,\Qr}\ot_V\Ga(\XX,\OO_{\XX}(u))^*\trig\Dcm_{\XX,\Qr}(-u)$. En utilisant la
proposition  ~\ref{subsubsection-dcm_gener}, on trouve ainsi l'énoncé
\begin{prop}\label{subsubsection-resol}
Soit $\NN$ un $\Dcm_{\XX,\Qr}$-module cohérent (resp. un
$\Ddag_{\XX,\Qr}$-module cohérent), alors 
il existe une résolution à deux termes $\Dcm_{\XX,\Qr}$-linéaire (resp. 
$\Ddag_{\XX,\Qr}$-linéaire) du type suivant 
 $$\left(\Dcm_{\XX}\right)^b\rig\left(\Dcm_{\XX}\right)^a\rig \NN\rig 0 \quad 
\left(resp. \, 
\left(\Ddag_{\XX,\Qr}\right)^b\rig\left(\Ddag_{\XX,\Qr}\right)^a\rig \NN\rig 0\right).$$
\end{prop} 
A partir de maintenant, $\DD$ désigne l'un des faisceaux $\Dcm_{\XX,\Qr}$ ou 
$\Ddag_{\XX,\Qr}$, $D=\Ga(\XX,\DD).$ Le corollaire suivant achève la démonstration 
du théorème. 
\begin{cor} \label{subsubsection-eq_cat}Les foncteurs $\Ga(\XX,.)$ et $\DD\ot_D .$ sont quasi-inverses et 
induisent une équivalence de catégories entre la catégorie des $D$-modules à gauche de présentation 
finie et la catégorie des $\DD$-modules à gauche cohérents. 
\end{cor} 
Soit $M$ un $D$-module de présentation finie et 
$$D^a\rig D^b\rig M\rig 0$$ une présentation de $M$. En tensorisant cette présentation 
par $\DD$, on trouve une présentation 
$$\DD^b\rig \DD^a\rig \DD\ot_D M\rig 0.$$
En particulier, le module $\DD\ot_D M$ est cohérent comme $\DD$-module à gauche. Par
acyclicité du foncteur $\Ga(\XX,.) $ pour les $\DD$-modules cohérents, on trouve un
diagramme dont les deux carrés sont commutatifs
$$\xymatrix{\Ga(\XX,\DD^b)\ar@{->}[r]&\Ga(\XX,\DD^a)\ar@{->}[r] & \Ga(\XX,\DD\ot_D M)
\ar@{->}[r] & 0 \\
 D^a \ar@{->}[u]^{\wr}\ar@{->}[r]& D^b \ar@{->}[u]^{\wr}\ar@{->}[r] & M\ar@{->}[u]^{g_M}\ar@{->}[r] & 0.}$$
Cela nous indique que la flèche $g_M$ est un isomorphisme. 

\vspace{+3mm}
Partons maintenant d'un $\DD$-module cohérent $\MM$. D'après la proposition précédente 
~\ref{subsubsection-resol}, il existe une résolution 
$$\DD^b \rig \DD^a \rig \MM \rig 0.$$ Comme le foncteur 
$\Ga$ est exact, on peut procéder comme précédemment pour voir 
que la flèche canonique $\DD\ot_D \Ga(\XX,\MM)\rig \MM$ est un isomorphisme.

Une première application de ce résultat est que l'on peut donner des propriétés de 
finitude des algèbres de sections globales des faisceaux d'opérateurs différentiels 
$\Dcm_{\XX,\Qr}$ et $\Ddag_{\XX,\Qr}$.
\subsection{Structure des algèbres de sections globales}
On garde les notations de~\ref{subsubsection-eq_cat}, en supposant toujours 
que $X$ vérifie (H) et que $X_K$ est $\DD_{X_K}$-affine. On a la proposition.
\begin{prop} \label{subsubsection-platitude} Le faisceau $\DD$ est un faisceau de $D$-modules plats à gauche.
\end{prop} 
\begin{proof}
Soit $I$ un idéal à gauche de présentation finie. Partons d'une suite exacte de 
$D$-modules à gauche de présentation finie
$$0\rig I \rig D \rig D/I \rig 0,.$$
On en déduit un complexe exact de $\DD$-modules à gauche cohérents, 
puisque le faisceau $\DD$ est cohérent
$$0 \rig Tor^1_D(\DD,D/I)\rig \DD\ot_D I \rig \DD \rig \DD/\DD I \rig 0.$$
On peut donc appliquer le foncteur exact $\Ga$, ce qui donne un diagramme dont tous les 
carrés sont commutatifs
$$\xymatrix{0\ar@{->}[r]&\Ga(\XX,Tor^1_D(\DD,D/I))\ar@{->}[r]&\Ga(\XX,\DD\ot_D I)\ar@{->}[r]
&\Ga(\XX,\DD)\ar@{->}[r]&\Ga(\XX,\DD/\DD I)\ar@{->}[r]&0 \\
0 \ar@{->}[r] & 0 \ar@{->}[u]\ar@{->}[r] & I \ar@{->}[r]\ar@{->}[u]^{\wr} &
D \ar@{->}[r]\ar@{->}[u]^{\wr}& D/I \ar@{->}[r]\ar@{->}[u]^{\wr}& 0 ,}$$
et dont les 3 dernières flèches verticales sont des isomorphismes en vertu de ~\ref{subsubsection-eq_cat}.
On en déduit que $\Ga(\XX,Tor^1_D(\DD,D/I))=0$ et, donc que 
$Tor^1_D(\DD,D/I))=0$ puisque ce module est cohérent. Supposons maintenant que 
$I$ est un idéal à gauche de type fini, alors il existe des idéaux à gauche de 
présentation finie, $\{I_i\}_{i\in \Omega}$ formant un système inductif et tel que 
$$I=\varinjlim_i I_i.$$ Comme le foncteur $Tor^1_D(\DD,.)$ commute à la limite inductive 
et que $$D/I=\varinjlim_i D/I_i,$$ on trouve
$$Tor^1_D(\DD,D/I)=\varinjlim_i Tor^1_D(\DD,D/I_i)=0,$$  
ce qui donne l'énoncé.
\end{proof}
A partir des propriétés de finitude des faisceaux $\Dcm_{\XX,\Qr}$ et $\Ddag_{\XX,\Qr}$,
on a l'énoncé suivant pour les schémas $\XX$ vérifiant les hypothèses du début de la
section.
\begin{thm} \be\item[(i)]La $V$-algèbre $\Ga(\XX,\Dcm_{\XX,\Qr})$ est une $V$-algèbre
complète noetherienne à gauche.
\item[(ii)]La $V$-algèbre $\Ga(\XX,\Ddag_{\XX,\Qr})$ est une $V$-algèbre faiblement complète 
cohérente à gauche. \ee
\end{thm} 
\begin{proof} Compte tenu de~\ref{subsubsection-coh_comp}, $\Ga(X,\Dcm_{\XX})$ est obtenu 
comme complété de $\Ga(X,\Dm_X)$, de sorte que $\Ga(X,\Dcm_{\XX,\Qr})$ est une 
$K$-algèbre de Banach. L'algèbre $\Ga(\XX,\Ddag_{\XX})$ est limite inductive (sur $m$) 
des algèbres $\Ga(X,\Dcm_{\XX})$ et est donc faiblement complète, tout comme 
$\Ga(\XX,\Ddag_{\XX,\Qr})$. Montrons maintenant la noetherianité de 
$\Ga(X,\Dcm_{\XX,\Qr})$. Notons comme précédemment $D=\Ga(X,\Dcm_{\XX,\Qr})$ et 
$\DD=\Dcm_{\XX,\Qr}$. Soit $I$ un idéal de $D$ et $\{I_i\}_{i\in \Omega}$ un système
inductif d'idéaux de présentation finie de $D$ tels que $$I=\varinjlim_i I_i.$$
introduisons $\II_i=\DD\ot_{D}I_i\subset \DD$, qui forment une suite croissante 
d'idéaux d'après l'énoncé de platitude ~\ref{subsubsection-platitude} précédent. 
Comme le faisceau $\DD$ est à sections noetheriennes sur les affines, il existe 
$i_0\in\Omega$ tel que $\II_{i_0}=\II_{j}$ pour tout $j\geq i_0$. Ce qui donne, en 
appliquant de nouveau~\ref{subsubsection-eq_cat} $I_{i_0}=I_{j}$ pour tout $j\geq
i_0$, et donc $I=I_0$ est de présentation finie et en particulier de type fini. En ce qui 
concerne (ii), considérons un idéal $I$ de type fini de $D=\Ga(X,\Ddag_{\XX,\Qr})$, 
et $J$ le noyau d'une surjection $D^a\rig I$. Ecrivons 
$$J=\varinjlim_i J_i,$$ où $\{J_i\}_{i\in \Omega}$ est un système inductif de 
$D$-modules à gauche de type fini. En tensorisant par $\DD$, en notant 
$\JJ_i=\DD\ot_D J_i$ (resp. $\JJ=\DD\ot_D J$), on trouve une
suite exacte de faisceaux de $\DD$-modules à gauche
$$ 0\rig \JJ \rig \DD^a \rig \DD\ot_D I \rig 0.$$
Comme le faisceau $\DD$ est cohérent sur les ouverts affines, il existe 
$i_0\in\Omega$ tel que $\forall j\geq i_0$, $\JJ_j=\JJ_{i_0}$. On conclut 
par les mêmes arguments que précédemment que $J=J_{i_0}$, de sorte que 
$I$ est de présentation finie.
\end{proof}

{\bf Remarque}: cet énoncé permet de voir que si $I$ est un $\Dcm_{\XX,\Qr}$-module de 
type fini, alors $\Ga(\XX,\Dcm_{\XX,\Qr}\ot_D I)\simeq I$ dans~\ref{subsubsection-eq_cat}.
%
%
\section{Le théorème de Beilinson-Bernstein arithmétique}
Tout est désormais en place pour le théorème principal de cet article. Reprenons les
notations de ~\ref{section-rappels}. Le théorème de Kempf de~\ref{paragraph-acyc_dom}, ainsi que 
~\ref{subsubsection-TX} entraînent que $X$ vérifie l'hypothèse (H). Soit $\lam\in X(T)$ et 
$\rho$ la demi-somme des racines positives de $G$. Si $\lam+\rho$ est dominant, 
le poids correspondant de l'algèbre de Lie $d\lam+d\rho$ est dominant
(\ref{subsubsection-comp_syst_rac}). Soit $\LL(\lam)$ 
le faisceau inversible associé à $\lam$. On peut donc appliquer le théorème 
principal de ~\cite{BeBe} et $X$ est $\DD_{X,\Qr}(\LL(\lam))$-affine. Pour simplifier les
notations, nous noterons dans la suite $\Dcm_{\XX,\Qr}(\lam)=\Dcm_{\XX,\Qr}(\LL(\lam))$
(resp. $\Ddag_{\XX,\Qr}(\lam)=\Ddag_{\XX,\Qr}(\LL(\lam))$).
Appliquons ~\ref{subsection-thm-int}. On obtient:
\begin{surthm} Soit $\lam\in X(T)$ tel que $\lam+\rho$ est dominant et régulier, alors 
$\XX$ est $\Ddag_{\XX,\Qr}(\lam)$-affine (resp. $\Dcm_{\XX,\Qr}(\lam)$-affine pour tout
entier $m$). 
\end{surthm} 
On obtient de plus les résultats suivants sur les algèbres de sections globales
sous les hypothèses du théorème précédent.
\begin{surthm} \be\item[(i)]Pour tout entier $m$, la $V$-algèbre
$\Ga(\XX,\Dcm_{\XX,\Qr}(\lam))$ est une $V$-algèbre
complète noetherienne à gauche.
\item[(ii)]La $V$-algèbre $\Ga(\XX,\Ddag_{\XX,\Qr}(\lam))$ est une $V$-algèbre faiblement complète 
cohérente à gauche. \ee
\end{surthm} 


\begin{thebibliography}{{A.~}67}

\bibitem[{A.~}67]{EGA4}
{A.~Grothendieck and {J}.~{D}ieudonn\'e}.
\newblock \'{E}l\'ements de {G}\'eom\'etrie {A}lg\'ebrique, etude locale des
  sch\'emas et des morphismes de sch\'emas, 4e partie.
\newblock {\em Publ. Math. I.H.E.S.}, 32, 1967.

\bibitem[BB81]{BeBe}
A.~Beilinson and J.~Bernstein.
\newblock { Localisation de ${\cal G}$-modules}.
\newblock {\em { Comptes-rendus Acad. Sc.}}, 292, p.~15--18, 1981.

\bibitem[Ber96]{Be1}
P.~Berthelot.
\newblock { $\DD$-modules arithm\'etiques I. {O}p\'erateurs diff\'erentiels de
  niveau fini}.
\newblock {\em Ann.~{scient}.~{\'E}c.~Norm.~Sup.}, $4^ e$ s\'erie, t.~29,
  p.185--272, 1996.

\bibitem[BK80]{Bry-Kash}
Jean-Luc Brylinski and Masaki Kashiwara.
\newblock D\'emonstration de la conjecture de {K}azhdan-{L}usztig sur les
  modules de {V}erma.
\newblock {\em C. R. Acad. Sci. Paris S\'er. A-B}, 291(6):373--376, 1980.

\bibitem[BLR90]{Neron_models}
Siegfried Bosch, Werner L{\"u}tkebohmert, and Michel Raynaud.
\newblock {\em N\'eron models}, volume~21 of {\em Ergebnisse der Mathematik und
  ihrer Grenzgebiete (3) [Results in Mathematics and Related Areas (3)]}.
\newblock Springer-Verlag, Berlin, 1990.

\bibitem[BMR04]{Bezru1}
Roman Bezrukavnikov, Ivan Mirkovic, and Dmitriy Rumynin.
\newblock {Localization of modules for a semisimple Lie algebra in prime
  characteristic}.
\newblock {\em {arXiv:math/0205144v8}}, 2004.

\bibitem[Bou68]{bourbaki_lie_456}
N.~Bourbaki.
\newblock {\em \'{E}l\'ements de math\'ematique. {F}asc. {XXXIV}. {G}roupes et
  alg\`ebres de {L}ie. {C}hapitre {IV}: {G}roupes de {C}oxeter et syst\`emes de
  {T}its. {C}hapitre {V}: {G}roupes engendr\'es par des r\'eflexions.
  {C}hapitre {VI}: syst\`emes de racines}.
\newblock Actualit\'es Scientifiques et Industrielles, No. 1337. Hermann,
  Paris, 1968.

\bibitem[DG70]{Demazure_gr_alg}
Michel Demazure and Pierre Gabriel.
\newblock {\em Groupes alg\'ebriques. {T}ome {I}: {G}\'eom\'etrie alg\'ebrique,
  g\'en\'eralit\'es, groupes commutatifs}.
\newblock Masson \& Cie, \'Editeur, Paris, 1970.
\newblock Avec un appendice {\it Corps de classes local}\ par Michiel
  Hazewinkel.

\bibitem[Haa87]{Haa}
B.~Haastert.
\newblock {\it {\"{U}}ber {D}ifferentialoperatoren und D-{M}oduln in positiver
  {C}harakteristik}.
\newblock {\em Manuscripta Mathematica}, 58, p.~385--415, 1987.

\bibitem[Huy97]{Hu1}
C.~Huyghe.
\newblock { $\DD^{\dagger}$-affinit\'e de l'espace projectif, avec un appendice
  de {P}. {B}erthelot}.
\newblock {\em Compositio Mathematica}, 108, No.~3, p.~277--318, 1997.

\bibitem[Jan03]{Jantzen}
Jens~Carsten Jantzen.
\newblock {\em Representations of algebraic groups}, volume 107 of {\em
  Mathematical Surveys and Monographs}.
\newblock American Mathematical Society, Providence, RI, second edition, 2003.

\bibitem[Kas89]{Kashi-D_flag}
Masaki Kashiwara.
\newblock Representation theory and {$D$}-modules on flag varieties.
\newblock {\em Ast\'erisque}, 173-174:9, 55--109, 1989.
\newblock Orbites unipotentes et repr\'esentations, III.

\bibitem[KL79]{Kazhdan-Lusztig1}
David Kazhdan and George Lusztig.
\newblock Representations of {C}oxeter groups and {H}ecke algebras.
\newblock {\em Invent. Math.}, 53(2):165--184, 1979.

\bibitem[KL02]{Kash-Laur}
Masaki Kashiwara and Niels Lauritzen.
\newblock Local cohomology and {$ D$}-affinity in positive characteristic.
\newblock {\em C. R. Math. Acad. Sci. Paris}, 335(12):993--996, 2002.

\bibitem[MFK94]{GIT}
D.~Mumford, J.~Fogarty, and F.~Kirwan.
\newblock {\em Geometric invariant theory}, volume~34 of {\em Ergebnisse der
  Mathematik und ihrer Grenzgebiete (2) [Results in Mathematics and Related
  Areas (2)]}.
\newblock Springer-Verlag, Berlin, third edition, 1994.

\end{thebibliography}
%

\noindent Christine Noot-Huyghe \\
\noindent Institut de Recherche Mathématique Avancée \\
\noindent Université Louis Pasteur et CNRS \\
\noindent 7, rue René Descartes \\
\noindent 67084 STRASBOURG cedex FRANCE \\
\noindent m\'el huyghe@math.u-strasbg.fr, http://www-irma.u-strasbg.fr/\textasciitilde huyghe

\end{document}